\documentclass[conference]{IEEEtran}
\IEEEoverridecommandlockouts
\usepackage{cite}
\usepackage{amsmath,amssymb,amsfonts,mathtools}
\usepackage{algorithmic}
\usepackage{graphicx}
\usepackage{textcomp}
\usepackage{xcolor}
\usepackage{caption}
\usepackage{subcaption}
\def\BibTeX{{\rm B\kern-.05em{\sc i\kern-.025em b}\kern-.08em
    T\kern-.1667em\lower.7ex\hbox{E}\kern-.125emX}}
\begin{document}

\title{Joint Expansion Planning of Power and Water Distribution Networks\\

\thanks{This work is funded by the Advanced Grid Modeling project of the US Department of Energy. The paper has been approved for unlimited public release and assigned LA-UR \# 23-33081.}
}

\author{\IEEEauthorblockN{Sai Krishna K. Hari, Russell Bent}
\IEEEauthorblockA{{Los Alamos National Laboratory}}
Los Alamos, NM, United States 
\and
\IEEEauthorblockN{Byron Tasseff}
\IEEEauthorblockA{Qcells North America} 
San Francisco, CA, United States \\
\and
\IEEEauthorblockN{Ahmed S. Zamzam, Clayton Barrows}
\IEEEauthorblockA{National Renewable Energy Laboratory} 
Golden, CO, United States 
% \and
% \IEEEauthorblockN{4\textsuperscript{th} Russell Bent}
% \IEEEauthorblockA{\textit{Applied Mathematics \& Plasma Physics} \\
% \textit{Los Alamos National Laboratory}\\
% Los Alamos, United States \\
% rbent@lanl.gov}
% \and
% \IEEEauthorblockN{5\textsuperscript{th} {\color{blue}Anyone else?}}
% \and
}

\maketitle

\begin{abstract}
This research explores the joint expansion planning of power and water distribution networks, which exhibit interdependence at various levels. We specifically focus on the dependency arising from the power consumption of pumps and develop models to seamlessly integrate new components into existing networks. Subsequently, we formulate the joint expansion planning as a Mixed Integer Nonlinear Program (MINLP). Through the application of this MINLP to a small-scale test network, we demonstrate the advantages of combining expansion planning, including cost savings and reduced redundancy, in comparison to independently expanding power and water distribution networks.
\end{abstract}


\section{Introduction}

Power and water networks are essential systems that are tightly coupled. The increasing consumption of these resources, unforeseen disasters, and scheduled maintenance drive the constant need to expand and upgrade the existing infrastructure networks delivering these essential resources. Due to the interdependence of these systems at various levels, such as hydropower generation, power plant cooling, water pumping, distribution, and water treatment, several recent articles have studied the integrated operations of these systems \cite{pereira2016joint, zamzam2018optimal, oikonomou2018optimal}. 
This operational interdependence heightens the vulnerability of these interconnected networks to cascading effects caused by natural and human-made disasters \cite{johnson2013building}. The question we address in this work is, ``Can we leverage the operational dependence of these interconnected networks to enhance the efficiency of their expansion and ultimately improve their resilience and operational efficiency?'"

To address this question, we develop models for incorporating various components into power and water distribution networks. These components include battery storage, generators, photovoltaic (PV) units, water tanks, pumps, and pipes. Furthermore, we model the interdependency between power and water networks arising from the power consumption of water pumps. Additionally, we present models governing the flow of power and water through components typical to these networks. Subsequently, we formulate a Mixed Integer Nonlinear Program (MINLP) using the high-fidelity nonlinear models discussed earlier. This formulation plays a crucial role in evaluating the advantages of joint expansion planning. We apply this formulation to a test network, effectively demonstrating the benefits derived from the integrated approach to expansion planning.

\section{Joint Network \& Expansion Modeling}
\label{sec:modeling}
We model the power and water distribution networks using graphs, where the vertices represent (power) buses and (water) junctions, and the edges represent bus-connecting or junction-connecting components such as power lines, water pipes, pumps, etc. To model the time evolution of various flow properties through these networks, we discretize the time horizon of planning into a set of time points, $\mathcal{K} = \{1, 2, \dots, K\}$. We denote the time intervals between these chosen time points by the set $\mathcal{K}^{\prime} = \{1, 2, \dots, K^{\prime} := K - 1\}$.

\subsection{Power Distribution Network Modeling}
\label{sec:power-models}
In power networks, we model the set of all buses, $\mathcal{N}$, as vertices and the set of all branches (power lines), $ \mathcal{E}$, as edges. We denote the set of all photovoltaic generators, battery storage units, and the loads that exist in the network or are available for addition to the network by $\mathcal{G}$, $\mathcal{B}$, and $\mathcal{D}$ respectively, and model them as elements connected to buses. We denote the subset of generators and storage units that do not already exist in the network but are available for expansion by $\widetilde{\mathcal{G}} \subset \mathcal{G}$ and $\widetilde{\mathcal{B}} \subset \mathcal{B}$.

% A power distribution network is represented by a graph $(\mathcal{N}, \mathcal{E})$, where $\mathcal{N}$ is the set of buses, $\mathcal{E}$ is the set of branches (lines).
% The sets of generators such as photovoltaic (PV) units, loads, and storage components are denoted by $\mathcal{G}$, $\mathcal{D}$, and $\mathcal{B}$, respectively.
% We let the subset of these components attached to $i \in \mathcal{N}$ be denoted by $\mathcal{G}_{i}$, $\mathcal{L}_{i}$, and $\mathcal{B}_{i}$.
% Furthermore, the set of conductors is denoted by $\mathcal{C}$ and is assumed to be equal to $\{a, b, c\}$.
% Additionally, there exists a subset of \emph{expansion} storage components and PV systems, which is indicated using a tilde.
% That is, $\widetilde{\mathcal{G}} \subset \mathcal{G}$ and $\widetilde{\mathcal{B}} \subset \mathcal{B}$ to denote the set of PV units and storage components \emph{that are to be considered for expansion}. 
% We next define the decision variables and constraints required to model an AC power distribution network's steady-state operations.

\subsubsection*{Buses}
% Each bus $i \in \mathcal{N}$ is associated with the complex voltage variables $V_{i, c}^{k} \in \mathbb{C}$, $c \in \mathcal{C}$, $k \in \mathcal{K}^{\prime}$.

We denote the voltage at each bus $i\in\mathcal{N}$ at time interval $k \in \mathcal{K}^{\prime}$ by $V_i^k := [V_{i,a}^k, V_{i,b}^k, V_{i,c}^k]^T$, and its upper and lower bounds by $\overline{v}_{i}^{k} := [\overline{v}_{i, a}^{k}, \overline{v}_{i, b}^{k}, \overline{v}_{i, c}^{k}]$ and $\underline{v}_{i}^{k} := [\underline{v}_{i, a}^{k}, \underline{v}_{i, b}^{k}, \underline{v}_{i, c}^{k}]$ respectively. Then, the operational limits on the voltage magnitude can be expressed as  
% The voltage at each bus $i\in\mathcal{N}$ at time interval $k \in \mathcal{K}^{\prime}$ is denoted by $V_i^k := [V_{i,a}^k, V_{i,b}^k, V_{i,c}^k]^T$. 
% The voltage magnitude bounds are represented as
% Operational limits require that each voltage magnitude resides between predefined lower and upper bounds, denoted by $\underline{v}_{i, c}^{k}$ and $\overline{v}_{i, c}^{k}$, respectively, i.e.,
\begin{equation}
    \underline{v}_{i}^{k} \leq \lvert V_{i}^{k} \rvert \leq \overline{v}_{i}^{k}, \, \forall i \in \mathcal{N}, \, \forall k \in \mathcal{K}^{\prime} \label{eqn:voltage-bounds}.
\end{equation}
% where $\underline{v}_{i}^{k} := [\underline{v}_{i, a}^{k}, \underline{v}_{i, b}^{k}, \underline{v}_{i, c}^{k}]$, and $\overline{v}_{i}^{k} := [\overline{v}_{i, a}^{k}, \overline{v}_{i, b}^{k}, \overline{v}_{i, c}^{k}]$ for conciseness. 

Let us denote the set of all edges connected to bus $i\in\mathcal{N}$ by $\mathcal{E}_i$, and the sets of generators, buses, and demand points connected to the bus by $\mathcal{G}_i$, $\mathcal{B}_i$, and $\mathcal{D}_i$ respectively. Besides, at the $k^{th}$ time interval, let the power injection provided by PV/generator $l \in \mathcal{G}_i$ be $S_{l}^{g, k}$ and the power withdrawal by energy storage $m \in \mathcal{B}_i$ and demand $n \in \mathcal{D}_i$ be $S_{m}^{b, k}$ and $S_{n}^{d, k}$ respectively, where $S_{l}^{g, k}, S_{m}^{b, k},  S_{n}^{d, k} \in \mathbb{C}^{3}$.
% During the $k$-th time interval at bus $i\in\mathcal{N}$, we denote the power injection provided by PV generator $l \in \mathcal{G}_i$ by  $S_{l}^{g, k}$ and power withdrawal by energy storage $m \in \mathcal{B}_i$ and demand $n \in \mathcal{D}_i$ by $S_{m}^{b, k}$, and $S_{l}^{g, k}, S_{n}^{d, k}$ respectively, where $S_{l}^{g, k}, S_{m}^{b, k},  S_{n}^{d, k} \in \mathbb{C}^{3}$.
% In addition, we use $ S_{i}^{d, k} \in \mathbb{C}^{3}$ to denote the fixed load demand at bus $i\in\mathcal{N}$ during the $k$-th time interval.
Then, following the Kirchoff's current law, the power flow balance at the buses can be expressed as 
% Then, the power flow balance at the buses can be expressed following the Kirchoff's current law as 
% % \subsubsubsection{Kirchoff's Current Law}
% Including bus-connected components, the nodal power balance equations can be expressed as
\begin{equation}
    \begin{aligned}
        &\sum_{\mathclap{\ell \in \mathcal{G}_{i}}} S^{g,k}_{\ell} -\sum_{\mathclap{\ell \in \mathcal{B}_{i}}} S^{b, k}_{\ell} - \sum_{\mathclap{\ell \in \mathcal{D}_{i}}} S^{d, k}_{\ell}\\
        &= \sum_{\mathclap{(i, j) \in \mathcal{E}_{i}}} \textnormal{diag}(S_{ij}^{k}), \, \forall i \in \mathcal{N}, \, \forall k \in \mathcal{K}^{\prime} \label{eqn:ac-kcl}.
    \end{aligned}
\end{equation}
% where $\mathcal{E}_{i}$ and $\mathcal{E}_{i}^{R}$ denote sets of branches adjacent to $i$.

\subsubsection*{Branches}
% Branches (lines) transport current throughout a power network.
% In a steady state, each branch satisfies Ohm's law, which relates the current $I_{ij}^{k} \in \mathbb{C}^{\mathcal{C} \times \mathcal{C}}$ through a branch to the voltages $V_{i}^{k}$ and $V_{j}^{k}$ at the buses connecting that branch.
% Ohm's law can instead be restated in terms of power, $S_{ij}^{k} \in \mathbb{C}^{\mathcal{C} \times \mathcal{C}}$, as

We denote the power flowing through a line  $(i,j)\in\mathcal{E}$ at time interval $k \in \mathcal{K}^{\prime}$ by $S_{ij}^{k}$. This power is related to the voltages at the buses connected by the line by Ohm's law, which can be expressed as 
% , and it is a function of voltages at the buses connected by the line. This relation is defined by the Ohm's law and can be expresed as 

% We denote the power flowing through line $(i,j)\in\mathcal{E}$ is denoted by $S_{ij}^{k}$ and express its relation to the voltage magnitidues 

% The power flowing on line $(i,j)\in\mathcal{E}$ is denoted by $S_{ij}^{k}$ and can be expressed as 
\begin{align}
    \begin{aligned}
        S_{ij}^{k} = \, & V_{i}^{k} (V_{i}^{k})^{H} (Y_{ij} + Y_{ij}^{c})^{H} \\
       -& V_{i}^{k} (V_{j}^{k})^{H} Y_{ij}^{H}, \, \forall (i, j) \in \mathcal{E}, \, \forall k \in \mathcal{K}^{\prime}
    \end{aligned} %\\
    % \begin{aligned}
    %     S_{ji}^{k} = \, & V_{j}^{k} (V_{j}^{k})^{H} (Y_{ij} + Y_{ji}^{c})^{H} \\
    %     - (& V_{i}^{k})^{H} V_{j}^{k} Y_{ij}^{H}, \, \forall (i, j) \in \mathcal{E}, \, \forall k \in \mathcal{K}^{\prime},
    % \end{aligned}
\end{align}
\label{eqn:ohms-law}%
where $Y_{ij} \in \mathbb{C}^{\mathcal{C} \times \mathcal{C}}$ is the admittance of the line; and $Y_{ij}^{c}, Y_{ji}^{c} \in \mathbb{C}^{\mathcal{C} \times \mathcal{C}}$ is line charging. 
The power flowing through the same line in the opposite direction can be expressed similarly and is omitted from discussion for conciseness.
% We can similarly express the power flowing on the same line in the opposite direction.
%For expansion branches $(i, j) \in \widetilde{\mathcal{E}}$, an indicator variable $x_{ij} \in \{0, 1\}$ is introduced to denote the expansion status of that component.
%That is, $x_{ij} = 1$ indicates the presence of that component in the expanded network and $x_{ij} = 0$ indicates its exclusion.
%For expansion branches, Ohm's law constraints are
%\begin{subequations}
%\begin{align}
%    \begin{aligned}
%        S_{ij}^{k} = \, & (\mathbf{e} x_{ij}) \circ [V_{i}^{k} (V_{i}^{k})^{H} (Y_{ij} + Y_{ij}^{c})^{H}] \\
%        -& (\mathbf{e} x_{ij}) \circ [V_{i}^{k} (V_{j}^{k})^{H} Y_{ij}^{H}], \, \forall (i, j) \in \widetilde{\mathcal{E}}, \, \forall k \in \mathcal{K}^{\prime}
%    \end{aligned} \\
%    \begin{aligned}
%        S_{ji}^{k} = \, & (\mathbf{e} x_{ij}) \circ [V_{j}^{k} (V_{j}^{k})^{H} (Y_{ij} + Y_{ji}^{c})^{H}] \\
%        - & (\mathbf{e} x_{ij}) \circ [(V_{i}^{k})^{H} V_{j}^{k} Y_{ij}^{H}], \, \forall (i, j) \in \widetilde{\mathcal{E}}, \, \forall k \in \mathcal{K}^{\prime},
%    \end{aligned}
%\end{align}
%\label{eqn:ohms-law-expansion}%
%\end{subequations}
%where $\mathbf{e} := \mathds{1}^{\mathcal{C} \times \mathcal{C}}$ and $\circ$ denotes the element-wise product.
%Note that the power flow along a branch is zero when $x_{ij} = 0$.
% Thermal limits are also imposed to prevent lines from sagging and protection devices from becoming active.

Thermal limits constraining the apparent power flow on the lines is given by
\begin{equation}
    \left\lvert \textnormal{diag}(S_{ij}^{k}) \right\rvert \leq \overline{s}_{ij}^{k}, \, \forall (i, j) \in \mathcal{E}, \, \forall k \in \mathcal{K}^{\prime} \label{eqn:thermal-limits-branch},
\end{equation}
where $\overline{s}_{ij}^{k} :=  (\overline{s}_{ij, a}^{k}, \overline{s}_{ij, b}^{k}, \overline{s}_{ij, c}^{k}) \in \mathbb{R}^{3}_{+}$ denote upper bounds.

% Finally, voltage phase angle differences are limited by predefined lower and upper bounds, $\underline{\theta}^{\Delta k}_{ij, c} \in \mathbb{R}$ and $\overline{\theta}^{\Delta k}_{ij, c} \in \mathbb{R}$, respectively, for all $(i, j) \in \mathcal{E}$, $c \in \mathcal{C}$ via the bound constraints
% \begin{equation}
%     \begin{gathered}
%         \underline{\theta}^{\Delta k}_{ij, c} \leq \angle \left(V_{i, c}^{k} V^{k *}_{j, c}\right) \leq \overline{\theta}^{\Delta k}_{ij, c}, \\
%         \forall (i, j) \in \mathcal{E}, \, \forall c \in \mathcal{C}, \, \forall k \in \mathcal{K}^{\prime}.
%     \end{gathered}
% \label{eqn:ac-phase-angle-difference-limits}%
% \end{equation}
%For expansion branches, these constraints are written as
%\begin{equation}
%    \begin{gathered}
%        -2 \pi (1 - x_{ij}) + \underline{\theta}^{\Delta k}_{ij, c} x_{ij} \leq \angle \left(V_{i, c}^{k} V^{k *}_{j, c}\right) \\
%        \leq \overline{\theta}^{\Delta k}_{ij, c} x_{ij} + 2 \pi (1 - x_{ij}), \,
%        \forall (i, j) \in \widetilde{\mathcal{E}}, \, \forall c \in \mathcal{C}, \, \forall k \in \mathcal{K}^{\prime}.
%    \end{gathered}
%\label{eqn:ac-phase-angle-difference-limits-expansion}%
%\end{equation}

% \subsection{Bus-connected Components}
% Bus-connected components comprise generators, loads, and shunts.
% Each is considered to be attached to a bus $i \in \mathcal{N}$.

\subsubsection*{Generators}
Suppose $\underline{S}_{\ell}^{g, k} \in \mathbb{C}$ and $\overline{S}_{\ell}^{g, k} \in \mathbb{C}$ denote the lower and upper bounds on complex power generation provided by generators or PV units. Then, the complex power supplied by these components is constrained as follows
% Generators or PV units inject power into the system, and we use $S_{\ell}^{g, k}$ to denote the complex power supplied by the unit $\ell \in \mathcal{G}$ during time interval $k \in \mathcal{K}^{\prime}$.
% Similarly, $\underline{S}_{\ell}^{g, k} \in \mathbb{C}$ and $\overline{S}_{\ell}^{g, k} \in \mathbb{C}$ denote the lower and upper bounds on complex power generation, respectively.
% This implies
\begin{equation}
    \underline{S}_{\ell}^{g, k} z_{\ell}^{g,k} \leq S_{\ell}^{g, k} \leq \overline{S}_{\ell}^{g, k} z_{\ell}^{g,k}, \, \forall \ell \in \mathcal{G}, \, \forall k \in \mathcal{K}^{\prime} \label{eqn:gen-bounds},
\end{equation}
where, with a slight abuse of notation, we use the vector inequalities here to bound the active and reactive injections at each phase separately. Here, $z_{\ell}^{g, k}$ is binary variable, which indicates the one/off (connected/disconnected) status of the generator when it is one and zero respectively. 

To model the expansion of generator $\ell\in\tilde{\mathcal{G}}$, we introduce a binary variable, $x_{\ell}^{g}$, and bound it by the on/off variable as shown below.
% We model the addition (expansion) of generator $\ell\in\tilde{\mathcal{G}}$ by introducing a binary variable, $x_{\ell}^{g}$, which takes the value one if the generator is added and zero if the generator is not added to the network. Then, the expansion of the generator can be captured by bounding the generator status variable, $z_{\ell}^{g, k}$, by the expansion variable, $x_{\ell}^{g}$, as shown below.
\begin{equation}
    \label{eqn:gen-expansion}
    z_{\ell}^{g,k} \leq x_{\ell}^g.\ \ \forall \ell\in\tilde{\mathcal{G}},\ k\in \mathcal{K}^\prime
\end{equation}
Here, the addition of the generator is indicated by $x_{\ell}^{g} = 1$. 
% Loads are points in the network where power is consumed.
% Herein, $S_{\ell}^{d, k} := (S_{\ell, a}^{d, k}, S_{\ell, b}^{d, k}, S_{\ell, c}^{d, k}) \in \mathbb{C}^{3}$ denote typically fixed but possibly variable demands at load $\ell \in \mathcal{B}$.

% \subsubsection{Shunts}
% Shunts are points in the network where current is measured or where a current-dependent voltage drop is desired.
% Herein, $Y_{\ell}^{s} \in \mathbb{C}^{\mathcal{C} \times \mathcal{C}}$ is the admittance of shunt $\ell \in \mathcal{H}$.

\subsubsection*{Storage}
% Storage components are points in the network where energy is stored and discharged over time.
We model each storage unit $\ell\in\mathcal{B}$ using a simplified version of the model from \cite{geth2020flexible}.
% Herein, $S_{\ell}^{b, k}$ denotes the power injection of the storage unit $\ell \in \mathcal{B}$ over time step $k \in \mathcal{K}^{\prime}$.
% Each power is partitioned into its active and reactive components as
% \begin{equation}
%     S_{\ell, c}^{u, k} = P_{\ell, c}^{u, k} + \mathbf{i} Q_{\ell, c}^{u, k}, \, \forall \ell \in \mathcal{U}, \, \forall c \in \mathcal{C}, \, \forall k \in \mathcal{K}^{\prime} \label{equation:power-storage-1},
% \end{equation}
We use $P_{\ell}^{b, k}$ and $Q_{l}^{b, k}$ to denote the active and reactive power injections into the storage respectively.
% The relationship between each power, current, and voltage magnitude is
% \begin{equation}
%     \lvert S_{\ell, c}^{u, k} \rvert^{2} = \lvert V_{i, c}^{k} \rvert^{2} \lvert I_{\ell, c}^{k} \rvert^{2}, \, \forall i \in \mathcal{N}, \, \forall \ell \in \mathcal{U}_{i}, \, \forall c \in \mathcal{C}, \, \forall k \in \mathcal{K}^{\prime} \label{equation:power-storage-2}.
% \end{equation}
Next, we introduce the variables $z_{\ell}^{b,k} \in \{0, 1\}$, $\ell \in \mathcal{B}$, $k \in \mathcal{K}^{\prime}$, to model the status of each energy storage at each time step, where zero and one indicate inactive and active statuses, respectively.
Then, bounds on the apparent and reactive power injections can be expressed using the following inequalities
% The apparent power injections are restricted by these variables via
% \begin{subequations}
    \begin{align}
        \lvert S_{\ell}^{b, k} \rvert &\leq z_{\ell}^{b,k} \overline{s}_{\ell}^{b}, \, \forall \ell \in \mathcal{B}, \, \forall k \in \mathcal{K}^{\prime},
        % \lvert I_{\ell, c}^{u, k} \rvert &\leq z_{\ell}^{k} \overline{i}_{\ell, c}^{u}, \, \forall \ell \in \mathcal{U}, \, \forall c \in \mathcal{C}, \, \forall k \in \mathcal{K}^{\prime},
    \end{align}
\label{equation:power-storage-3}%
% \end{subequations}
% {\color{red}where $\overline{s}_{\ell, c}^{u}$ is the apparent power rating of the battery inverter}. The reactive power variables $Q_{\ell}^{b, k}$ are bounded as per
\begin{equation}
    -\overline{s}_{\ell}^{b} z_{\ell}^{b, k} \leq Q_{\ell}^{b, k} \leq \overline{s}_{\ell}^{b} z_{\ell}^{b, k}, \, \forall \ell \in \mathcal{B}, \, \forall k \in \mathcal{K}^{\prime} \label{equation:power-storage-4},
\end{equation}
where $\overline{s}_{\ell}^{b}$ is the apparent power rating of the battery inverter.
% Next, let $P_{\ell}^{b k} = P_{\ell}^{u, f, k} - P_{\ell}^{u, g, k}$ denote the power from the storage subcomponent of $u \in \mathcal{U}$, where $P_{\ell}^{u, f, k}, \, P_{\ell}^{u, g, k} \geq 0$ denote storage charging and discharging, respectively.
% Since storage charging and discharging are typically complementary,
% \begin{equation}
%     P_{\ell}^{u, f, k} P_{\ell}^{u, g, k} = 0, \, \forall \ell \in \mathcal{U}, \, \forall k \in \mathcal{K}^{\prime} \label{equation:power-storage-5}.
% \end{equation}
Next, we model the energy content of the storage buffer using the following equation, 
\begin{gather}
    E_{\ell}^{k + 1} = E_{\ell}^{k} + \Delta t^{k} \left(\eta_{\ell}^{c} P_{\ell}^{bc,k + 1} - \frac{1}{\eta_{\ell}^{d}}P_{\ell}^{bd, k + 1} \right), 
    % \forall \ell \in \mathcal{U}, \, \forall k \in \mathcal{K}^{\prime} 
    \label{equation:power-storage-6}
\end{gather}
where $ P_{\ell}^{bc,k + 1}$ and $P_{\ell}^{bd,k + 1}$ indicate the charging and discharging power, the difference and products of which are $P_{\ell}^{b,k + 1}$ and 0 respectively.
Here, $E_{\ell}^{k}$ is the variable energy stored at time index $k \in \mathcal{K}$, $\Delta t^{k}$ is the time step between the discrete time points $k \in \mathcal{K}$ and $k + 1 \in \mathcal{K}$, $0 \leq \eta_{\ell}^{c}, \eta_{\ell}^{d} \leq 1$ are charging and discharging efficiencies.
To ensure stored energy is recovered by the end of the planning period, we also enforce the following constraints.
\begin{equation}
    E_{\ell}^{K} \geq E_{\ell}^{1}, \, \forall \ell \in \mathcal{B} \label{equation:power-storage-7}
\end{equation}

% One could also enforce fixed values of $E_{\ell}^{u, 1}$.
Additionally, we enforce the following operational constraints, which limit the charging and discharging rates, and the energy stored in the battery by specified upper bounds, $\overline{P}_{\ell}^{c}$, $\overline{P}_{\ell}^{d}$, and $\overline{E}_{\ell}$ respectively.
\begin{subequations}
    \begin{align}
        0 \leq P_{\ell}^{c} \leq \overline{P}_{\ell}^{c}z_{\ell}^{b, k}, \, &\forall \ell \in \mathcal{B}, \, \forall k \in \mathcal{K}^{\prime}, \\
        0 \leq P_{\ell}^{d} \leq \overline{P}_{\ell}^{d}z_{\ell}^{b, k}, \, &\forall \ell \in \mathcal{B}, \, \forall k \in \mathcal{K}^{\prime}, \\
        0 \leq E_{\ell}^{k} \leq \overline{E}_{\ell}, \, &\forall \ell \in \mathcal{B}, \, \forall k \in \mathcal{K}^{\prime}.
    \end{align}
\label{equation:power-storage-8}%
\end{subequations}
% where $\overline{P}_{\ell}^{c}$, $\overline{P}_{\ell}^{d}$, and $\overline{E}_{\ell}$ are upper bounds on the charging and discharging rate and the energy stored.

% Combining properties of the converter and storage subcomponent, the power balance constraints for all components are
% \begin{equation}
% \begin{gathered}
%     P_{\ell}^{u, s, k} + \sum_{c \in \mathcal{C}} S_{\ell, c}^{u, k} = \mathbf{i} Q_{\ell}^{u, n, k} + S_{\ell}^{e, k} \\
%     + \sum_{c \in \mathcal{C}} Z_{\ell, c}^{u} \lvert I_{\ell, c}^{u, k} \rvert^{2}, \,
%     \forall \ell \in \mathcal{U}, \, \forall c \in \mathcal{C}, \, \forall k \in \mathcal{K}^{\prime}.
% \end{gathered}
% \label{equation:power-storage-9}%
% \end{equation}
% Here, $S_{\ell}^{e, k} \in \mathbb{C}$ represents constant external power loss at $k \in \mathcal{K}^{\prime}$, and the last sum represents copper losses, where $Z_{\ell, c}^{u}$ is the internal impedence of the converter for conductor $c \in \mathcal{C}$.

Finally, we model the addition of expanded storage by introducing a binary variable $x^b_{\ell} \in \{0, 1\}$, $\ell \in \widetilde{\mathcal{B}}$, and constraining the on/off status variable as follows
% letting $x_{\ell} \in \{0, 1\}$, $\ell \in \widetilde{\mathcal{B}}$, denote the expansion decisions for each expansion storage component, the following constraints are used to model inclusion of expanded storage:
\begin{equation}
    z_{\ell}^{b,k} \leq x^b_{\ell}, \, \forall \ell \in \widetilde{\mathcal{B}}, \, \forall k \in \mathcal{K}^{\prime}.
\label{equation:power-storage-10}%
\end{equation}

\subsection{Water Distribution Network Modeling}
\label{sec:water-models}
% \subsection{Notation for Sets}

% Water distribution networks 
In water networks, we represent the junctions by vertices, and the elements connecting junctions, such as pipes and pumps, by edges. The sets of all pipes and pumps in the network are denoted by $\mathcal{A}$ and $\mathcal{P}$ respectively, and they constitute the set of all edges, $\mathcal{L}$. We model reservoirs, tanks, and demand points as elements connected to junctions and denote the sets of these elements by $\mathcal{R}$, $\mathcal{T}$, and $\mathcal{W}$ respectively. The set of all junctions is denoted by $\mathcal{J}$. For a compact representation of models, we denote the subsets of reservoirs, tanks, and demands attached to node $i \in \mathcal{J}$ by $\mathcal{R}_{i} \subset \mathcal{R}$, $\mathcal{T}_{i} \subset \mathcal{T}$, and $\mathcal{W}_{i} \subset \mathcal{W}$ respectively.
We denote the subsets of all the components considered for expansion by a tilde ($\sim$) over the set. For example, $\widetilde{\mathcal{L}}$ represents the set of all edges considered for expansion and satisfies $\widetilde{\mathcal{L}} \subset \mathcal{L}$.

\subsubsection*{Junctions / Nodes}
We denote the total hydraulic head (also referred to as ``head'') at a junction $i \in \mathcal{J}$ at time instance $k \in \mathcal{K}$ by $h_{i}^{k}$. These variables are subject to bounds on the minimum and maximum head, $\underline{h}_{i}^{k}$ and $\overline{h}_{i}^{k}$, which are often inferred from the network data. These bounds take the form of the following set of inequalities
% Nodal potentials are denoted by the variables $h_{i}^{k}$, where each variable represents the total hydraulic head (hereafter referred to as ``head'') in units of length at a time instance.
% The total hydraulic head  assimilates elevation and pressure heads at a node, while the velocity head is neglected.
% For each $i \in \mathcal{J}$, $k \in \mathcal{K}$, a minimum head $\underline{h}_{i}^{k}$, determined a priori, must first be satisfied.
% Additionally, an upper bound $\overline{h}_{i}^{k}$ can often be inferred from network data.
% This implies the head bounds
\begin{equation}
    \underline{h}_{i}^{k} \leq h_{i}^{k} \leq \overline{h}_{i}^{k}, \, \forall i \in \mathcal{J}, \; \forall k \in \mathcal{K} \label{equation:node-head-bounds}.
\end{equation}

\subsubsection*{Edges}
Every node-connecting component $(i, j) \in \mathcal{L}$ is associated with a variable, $q_{ij}^{k}$, which denotes the volumetric flow rate across that component.
Note that flow can assume both positive and negative values.
Assuming lower and upper bounds of $\underline{q}_{ij}^{k}$ and $\overline{q}_{ij}^{k}$, respectively, non-expansion components must first satisfy the bounds
\begin{equation}
    \underline{q}_{ij}^{k} \leq q_{ij}^{k} \leq \overline{q}_{ij}^{k}, \, \forall (i, j) \in \mathcal{L} \setminus \widetilde{\mathcal{L}}, \, \forall k \in \mathcal{K}^{\prime} \label{equation:non-expansion-flow-bounds}.
\end{equation}

For each expansion component $(i, j) \in \widetilde{\mathcal{L}}$, we introduce an indicator variable $x_{ij}^{w} \in \{0, 1\}$ to denote the expansion status of that component. 
% Here, $x_{ij} = 1$ indicates the presence of that component in the expanded network and $x_{ij} = 0$ indicates its exclusion.
Therefore, the conditional flow bounds can be expressed as
\begin{equation}
    \underline{q}_{ij}^{k} x_{ij} \leq q_{ij}^{k} \leq \overline{q}_{ij}^{k} x_{ij}, \, \forall (i, j) \in \widetilde{\mathcal{L}}, \, \forall k \in \mathcal{K}^{\prime} \label{equation:expansion-flow-bounds}.
\end{equation}

\subsubsection*{Pipes}
The flow of water along a pipe is induced by the head difference between the nodes connected by the pipe. The relation between the flow rate and the head loss is typically nonlinear, and here, we capture this nonlinear relation using the Hazen-Williams equation \cite{ormsbee2016darcy}. Then, the head loss along pipes that are existing in the network can be modeled using the following equation
\begin{equation}
    h_{i}^{k} - h_{j}^{k} = L_{ij} r_{ij} q_{ij}^{k} \lvert q_{ij}^{k} \rvert^{.852}, \; \forall (i, j) \in \mathcal{A} \setminus \widetilde{\mathcal{A}}, \, \forall k \in \mathcal{K}^{\prime}, \label{equation:non-expansion-pipe-head-loss}.
\end{equation}
where $L_{ij}$ denotes the pipe length and $r_{ij}$ denotes the pipe resistance per unit length.
% as $r_{ij} = \frac{10.7} {\kappa_{ij}^{1.852} D_{ij}^{4.8704}}$.
% $\tau_{ij}^{k}$ is the friction factor, $g$ is gravitational acceleration, and $D_{ij}$ is the diameter.}

For pipes considered for expansion, we use the binary variable $x_{ij}$ to indicate the expansion status. If the pipe is added to the network, the Hazen-Williams equation must apply and if the pipe is not added, the flow through it must be zero and the heads at either nodes must be decoupled. This can be modeled using the following equation
\begin{equation}
    x_{ij} (h_{i}^{k} - h_{j}^{k}) = L_{ij} r_{ij} q_{ij}^{k} \lvert q_{ij}^{k} \rvert^{0.852}, \; \forall (i, j) \in \widetilde{\mathcal{A}}, \, \forall k \in \mathcal{K}^{\prime} \label{equation:expansion-pipe-head-loss}.
\end{equation}

\subsubsection*{Pumps}
Here, we consider only fixed-speed pumps. We use the binary variable, $z_{ij}^{p,k}$, to indicate the on/off status of pump $(i, j) \in \mathcal{P}$ at time interval $k$. We assume that the pump permits only unidirectional flow and the flow is non-zero when the pump is active. We model this using the following constraint.
% provides a non-zero head increase from node $i$ to $j$ when active and zero flow when inactive.
% and permits a unidirectional flow.
% Each pump $(i, j) \in \mathcal{P}$ increases the head from node $i$ to $j$ when active and permits only unidirectional flow.
% Here, we consider only fixed-speed pumps, and indicate the on/off status of each pump at time step $k$ by the binary variable $z_{ij}^{p,k}$. Then, 
% where each pump is assumed to be either on or off.
% When the pump is off, there is zero flow along the pump, and heads at adjacent nodes are decoupled.
% When the pump is on, there is appreciably positive flow (greater than or equal to some fixed $\epsilon_{ij}$), and the head increase from $i$ to $j$ is described by a nonlinear function.
% The variable $z_{ij}^{p,k} \in \{0, 1\}$ indicates the on/off status of each pump. %, where $z_{ij}^{k} = 1$ if $q_{ij}^{k} \geq \epsilon_{ij}$ and $z_{ij}^{k} = 0$ if $q_{ij}^{k} < \epsilon_{ij}$.
% This is modelled as
\begin{equation}
    \underline{q}_{ij}^{k} = 0 \leq \epsilon_{ij} z_{ij}^{p,k} \leq q_{ij}^{k} \leq \overline{q}_{ij}^{k} z_{ij}^{p,k}, \, \forall (i, j) \in \mathcal{P}, \, \forall k \in \mathcal{K}^{\prime} \label{equation:pump-flow-bounds}.
\end{equation}

The pump is assumed to provide a non-zero head gain when active, and this gain is denoted by the variable $G_{ij}^{k} \geq 0$. Following \cite{ulanicki2008modeling}, we model the head gain provided by the pump as a strictly concave function of the following form
% Modeling a pump's head gain as a quadratic function, when active, is an established practice that is well-documented in the literature (e.g., \cite{ulanicki2008modeling}).
% We assume that each head gain is modeled via a strictly concave function of the form
\begin{equation}
    a_{ij} \left(q_{ij}^{k}\right)^{2} + b_{ij} q_{ij}^{k} + c_{ij} z_{ij}^{p,k} = G_{ij}^{k},  \, \forall (i, j) \in \mathcal{P}, \, \forall k \in \mathcal{K}^{\prime} \label{equation:pump-head-gain}.
\end{equation}
Here, $a_{ij} < 0$, and the head gain takes a minimum value of $c_{ij} > 0$ when the pump is active. When the pump is off, we use similar constraints to decouple the heads at the nodes connected by the pumps. We omit these constraints from the discussion in the interest of space.
% Note further that $q_{ij}^{k}$ and $c_{ij} z_{ij}^{k}$ are restricted to zero when $z_{ij}^{p,k}$ is zero.
% This ensures that, when a pump is off, the corresponding head gain $G_{ij}^{k}$ is also zero.
% To ensure the decoupling of hydraulic heads when a pump is off, the following constraints are appended:
% \begin{subequations}
% \begin{align}
%     & \begin{gathered}
%         h_{j}^{k} - h_{i}^{k} \leq G_{ij}^{k} + \left(1 - z_{ij}^{k}\right) (\overline{h}_{j}^{k} - \underline{h}_{i}^{k}), \\ \forall (i, j) \in \mathcal{P}, \, \forall k \in \mathcal{K}^{\prime}
%     \end{gathered} \\
%     & \begin{gathered}
%         h_{j}^{k} - h_{i}^{k} \geq G_{ij}^{k} + \left(1 - z_{ij}^{k}\right) (\underline{h}_{j}^{k} - \overline{h}_{i}^{k}), \\ \forall (i, j) \in \mathcal{P}, \, \forall k \in \mathcal{K}^{\prime}.
%     \end{gathered}
% \end{align}
% \label{equation:pump-pressure}%
% \end{subequations}
% Note that when $z_{ij}^{k} = 1$, the pump is on, and the head gain between nodes is $G_{ij}^{k}$.
% Otherwise, the heads are decoupled.
% {\color{blue}Similar constraints are also introduced to ensure heads are decoupled when pumps are off.}
For pumps considered for expansion, the binary variable $x_{ij}$ limits their activations as per
\begin{equation}
    z_{ij}^{p,k} \leq x_{ij}, \, \forall (i, j) \in \widetilde{\mathcal{P}}, \, \forall k \in \mathcal{K}^{\prime}.
\label{equation:pump-expansion}%
\end{equation}

% \subsection{Node-connected Components}
% Different from node-\emph{connecting} components, node-\emph{connected} components comprise demands, reservoirs, and tanks.
% Each is considered to be attached to a node $i \in \mathcal{J}$.

\subsubsection*{Demands \& Reservoirs}
Each demand point  $n \in \mathcal{W}$ in the network is associated with a variable $q^{k}_{n}$, which denotes the volumetric rate at which water is removed from the point at time $k \in \mathcal{K}^{\prime}$. On the other hand, reservoirs are assumed to be infinite sources of supply that can inject water into the network at a flow rate $q_{n}^{k} \geq 0$, $ \forall n \in \mathcal{R}$, $k \in \mathcal{K}^{\prime}$.

% The reservoir on the other hand, is assumed to be an infinite source of flow with zero pressure and constant elevation at a given time point.
% Demands are discrete points in the network attached to nodes where water is typically withdrawn.
% Each demand in the network is associated with a variable $\overline{q}^{k}_{n}$ that denotes the demanded flow, expressed as a volumetric flow rate, at demand $n \in \mathcal{W}$, time $k \in \mathcal{K}^{\prime}$.
% Note that, without loss of generality, demands are allowed to be negative.
% A negative demand indicates a point of injection of water into the network.
% For convenience, the variables $q_{n}^{k} \in \mathbb{R}$, $n \in \mathcal{D}$, $k \in \mathcal{K}^{\prime}$, are introduced to denote the amount of water consumed by each demand point.
% (Note that when demand is fixed, $q_{n}^{k} = \overline{q}^{k}_{n}$.)
% \subsubsection{Reservoirs}
% Reservoirs are discrete points in the network attached to nodes where water is supplied.
% Also, each reservoir is assumed to be an infinite source of flow with zero pressure and constant elevation at a point in time (i.e., $\underline{h}_{i}^{k} = \overline{h}_{i}^{k}$ is assumed at the attached node).
% The variables $q_{n}^{k} \geq 0$, $n \in \mathcal{R}$, $k \in \mathcal{K}^{\prime}$, denote the outflow of water from each reservoir at time $k$.

\subsubsection*{Tanks}
We model tanks as cylindrical water storage facilities with a constant cross-sectional area. For tank $n \in \mathcal{T}$, we denote the cross-sectional area by $A_n$ and model the volume of stored water at time point $k \in \mathcal{K}$ by  $ V_{n}^{k} $ and the rate of water withdrawal from the tank between time points $k$ and $k + 1 \in \mathcal{K}$ by $ q_{n}^{k} $.
Then, the constraints governing the flow of water through the tanks is given by \eqref{equation:tank-volume-integration} and \eqref{equation:tank-volume-expression}.
% Tanks are discrete points in the network attached to nodes that serve as means for storing and discharging water over time.
% In this study, all tanks are assumed to be cylindrical with a fixed diameter $D_{n}$, where $n \in \mathcal{T}$, which implies a cross-sectional area $A_{n} := \frac{\pi}{4} D_{n}^{2}$.
% The bottom of each tank is assumed to be located at or below the minimum elevation of the connecting node, $b_{i}^{k} \leq \underline{h}_{i}^{k}$, $i \in \mathcal{J}$, and the maximum elevation of water in the tank is assumed to be $\overline{h}_{i}^{k}$.
% The variables $q_{n}^{k}$, $n \in \mathcal{T}$, $k \in \mathcal{K}^{\prime}$, denote the outflow through each tank.
% With these variables, volumes are first defined as
\begin{equation}
    V_{n}^{k+1} = V_{n}^{k} - \Delta t^{k} q_{n}^{k}, \, \forall n \in \mathcal{T}, \, \forall k \in \mathcal{K} \label{equation:tank-volume-integration},
\end{equation}
\begin{equation}
    V_{n}^{k} := A_{n} (h_{i}^{k} - b_{i}^{k}), \, \forall i \in \mathcal{J}, \, \forall n \in \mathcal{T}_{i}, \, \forall k \in \mathcal{K} \label{equation:tank-volume-expression},
\end{equation}
where $b_i^{k}$ indicates the minimum elevation of the node connected to the tank.
We also require the water level in the tanks to recover to their initial level at least and model this using the following constraint
% We also assume that the tank volume at the end of the horizon is at least equal to its volume in the first time step. This is modeled as follows
% The Euler steps for integrating tank volumes are then given as
% One important criterion for an expansion is that, by the end of the scheduling horizon, the volume of water within each tank will be at least as large as its initial volume.
% This ensures that (i) the resultant optimal component schedule (e.g., of pump activations) can be repeatedly applied to subsequent planning periods with similar demand profiles and (ii) tank volumes at the end of the planning period are not small enough to discourage efficient operations in subsequent periods.
% The tank volume recovery implies the constraints
\begin{equation}
    V_{n}^{1} \leq V_{n}^{K}, \, \forall n \in \mathcal{T} \label{equation:tank-volume-recovery}.
\end{equation}
We model the expansion of tanks by introducing a small node-connecting element, such as a short pipe, between the tank and the junction to which it is associated. Then, the tank is added to the network iff the pipe is also added. The addition of tank $n \in \widetilde{\mathcal{T}}$ to the network is modeled by the binary variable $x_n$.

\subsubsection*{Flow Conservation}
Finally, the flow conservation at every junction is modeled as follows
% Satisfaction of flow conservation throughout the network requires flow balance constraints to be enforced at every node $i \in \mathcal{J}$ across all intervals $\mathcal{K}^{\prime}$, i.e.,
\begin{equation}
    \sum_{\mathclap{j:(j, i)\in\mathcal{L}}} q_{ji}^{k}  = \sum_{\mathclap{n \in \mathcal{W}_{i}}} q_{n}^{k} - \sum_{\mathclap{n \in \mathcal{R}_{i}}} q_{n}^{k} - \sum_{\mathclap{n \in \mathcal{T}_{i}}} q_{n}^{k}, \, \forall i \in \mathcal{J}, \, \forall k \in \mathcal{K}^{\prime} \label{equation:flow-conservation}.
\end{equation}

\subsection{Joint Network Modeling}
\label{sec:joint-modles}
The primary interdependency between power and water distribution systems is through the operation of pumps, which consume electricity.
The power consumption of pumps is a nonlinear function of their efficiency, water flow rate through them, and the head gain they offer
% flow rate, and the head gain offered 
\cite{ulanicki2008modeling}. 
However, in the case of fixed-speed pumps, exhibiting power curves akin to those outlined in \cite{ulanicki2008modeling}, power consumption can be effectively characterized by an affine function of the flow rate, as articulated in Eq. \eqref{equation:pump-power-linear}.

% The power consumption of pump $(i, j) \in \mathcal{P}$ at time $k \in \mathcal{K}^{\prime}$ is accurately defined \cite{ulanicki2008modeling} via
% \begin{equation}
%     P_{ij}^{k}(q_{ij}^{k}, G_{ij}^{k}) := \frac{\rho g G_{ij}^{k} q_{ij}^{k}}{\eta_{ij}(q_{ij}^{k})}, \, \forall (i, j) \in \mathcal{P}, \, \forall k \in \mathcal{K}^{\prime}.
% \end{equation}
% Here, $P_{ij}^{k}$ is the pump power consumed, $\rho$ is the density of water, $g$ is the acceleration due to gravity, and $\eta_{ij}(\cdot)$ is the pump's dimensionless efficiency, which is often a nonlinear function of flow and dictated by a predefined curve provided by the pump manufacturer.
% Note that each equation for a pump's power consumption is thus highly nonlinear.
% As such, some optimization studies have resorted to power curves for fixed speed pumps that are linear as a function of flow.
% {\color{blue} Any references for this?}
% This results in power functions parameterized as, for example,
\begin{equation}
    P_{ij}^{k}(q_{ij}^{k}, z_{ij}^{k}) := \alpha_{ij} q_{ij}^{k} + \mu_{ij} z_{ij}^{k}, \, \forall (i, j) \in \mathcal{P}, \, \forall k \in \mathcal{K}^{\prime},
    \label{equation:pump-power-linear}%
\end{equation}
where $\alpha_{ij}$ and $\mu_{ij}$ are fixed constants determined from the power curves. Here, $\mu_{ij}$ describes the non-zero power consumption of pumps when active.
% Pumps in the water network are naturally modeled as variable loads in the power network.

Let $\mathcal{I}$ denote the set of interdependencies that link each pump $(i, j) \in \mathcal{P}$ with a load $\ell \in \mathcal{B}$. Then, the constraints linking the power and water networks are modeled by spreading the active power consumption of pumps across that of all conductors of the power load as follows
\begin{equation}
    P_{ij}^{k}(q_{ij}^{k}, \cdot) = \sum_{c \in \mathcal{C}} \Re(S_{\ell, c}^{d, k}), \, \forall (i, j, \ell) \in \mathcal{I}, \, \forall k \in \mathcal{K}^{\prime} \label{equation:pump-as-load},
\end{equation}
where $\Re(S_{\ell, c}^{d, k})$ denotes the active power of conductor $c$ of load $\ell$.
% {\color{blue} Budget constraint}

\section{Mathematical Formulation for Joint Expansion Planning}
In this work, we consider expansion planning with the objective of maximizing the weighted sum of water demands and power demands (excluding pump loads) satisfied by the expanded network subject to a given budget on the cost of expansion. We use the parameter $0 \leq \lambda \leq 1$ to model the priority given to power demand over its counterpart. Then, the problem objective function takes the following form

% % Depending on the specific goal of expansion, one can choose from a wide variety of metrics ranging from maximizing the demand satisfied by the expanded networks to increasing the 
% such as increasing the demand met or improving the resilience of the networks to threats, the objective functions for the joint problem vary. For example, one can either minimize the weighted cost of expanding the power and water networks for given expected demands that must be satisfied or maximize the weighted power and water demands satisfied subject to a constraint on the cost of expansion. In this work, we choose the later. Thus, the objective function takes the form of 
% {\color{blue} Scaling. Same sets used for power and water demand points.}
\begin{equation}
    f(\cdot) = \sum_{k \in \mathcal{K}} (\lambda \sum_{l \in \mathcal{D
    } \setminus \mathcal{LI}}S^{d,k}_{l} + (1 - \lambda) \sum_{i \in \mathcal{W}}q^k_{i}) ,
    \label{eq:objective}
\end{equation}
where the set $\mathcal{LI}$ indicates the loads connected to pumps.

The constraint on the cost of expansion, which includes the cost of adding storage units, generators, water tanks, pipes, and pumps, can be expressed as
\begin{equation}
      \sum_{\mathclap{\ell \in \widetilde{\mathcal{B}} \cup \widetilde{\mathcal{G}} \cup \widetilde{\mathcal{T}}}} c_{\ell} x_{\ell} + \sum_{\mathclap{(i,j) \in \widetilde{\mathcal{A}} \cup \widetilde{\mathcal{P}}}} c_{(i,j)} x_{(i,j)} \leq \textbf{B},
      \label{ineq:budget}
\end{equation}
where $c_{(.)}$ indicates the cost of adding a component to the network, and \textbf{B} indicates the budget available for expansion.

% Two network expansion objectives have thus far been considered: minimization of power network costs, i.e., $C_{\mathcal{U} \cup \mathcal{G}}$, and minimization of water network costs, i.e., $C_{\mathcal{A} \cup \mathcal{S}} + C_{\mathcal{V} \cup \mathcal{W} \cup \mathcal{P}}$.
% In practice, one objective may take priority over the other.
% Here, we choose to model the tradeoff between power and water network costs via a nonnegative scaling parameter, $0 \leq \lambda \leq 1$, and where the objective function is defined as
% \begin{equation}
%     f(\cdot) = \lambda (C_{\mathcal{U} \cup \mathcal{G}}) + (1 - \lambda) (C_{\mathcal{A} \cup \mathcal{S}}+ C_{\mathcal{V} \cup \mathcal{W} \cup \mathcal{P}}).
% \end{equation}
% For example, when $\lambda = 1$, power network costs are prioritized, and when $\lambda = 0$, water network costs are prioritized.

% \subsection{Joint Optimal Network Expansion}
Then, the joint expansion problem can be formulated as
\begin{equation*}
    \begin{aligned}
        & \text{Maximize} & & \textnormal{Weighted sum of power \& water demands: } \eqref{eq:objective}\\
        & \text{subject to} & & \textnormal{Power network constraints}: \eqref{eqn:voltage-bounds} - \eqref{equation:power-storage-10}\\
        & & & \textnormal{Water network constraints}: \eqref{equation:node-head-bounds}\ - \eqref{equation:flow-conservation}, \eqref{equation:pump-power-linear}\\
        & & & \textnormal{Pump-as-load constraints}: \eqref{equation:pump-as-load}\\
        & & & \textnormal{Budget constraints}: \eqref{ineq:budget}.
    \end{aligned}
\end{equation*}

% Maximize \eqref{eq:objective}, subject to \eqref{equation:ncp-power} (power network constraints), \eqref{equation:mincp-f} (water network constraints), \eqref{equation:pump-as-load} (linking constraints modeling pumps as loads in power networks), and \eqref{budget} (budget constraints).
Due to the presence of nonlinear constraints and discrete variables, the formulation takes the form of an MINLP.
% \begin{equation}
% \begin{aligned}
%     % & {\color{red}\text{minimize}} & & \lambda (C_{\mathcal{U} \cup \mathcal{G}}) + (1 - \lambda) (C_{\mathcal{A} \cup \mathcal{S}}+ C_{\mathcal{V} \cup \mathcal{W} \cup \mathcal{P}}) \\
%     & \text{\color{blue}{maximize}} & & \lambda \text{(power demand)} + (1 - \lambda) \text{(water demand)}\\
%     & \text{subject to} & & \textnormal{Power network Constraints} ~ \eqref{equation:ncp-power} \\
%     & & & \textnormal{Water network Constraints} ~ \eqref{equation:mincp-f} \\
%     & & & \textnormal{Pump-as-load Constraints} ~ \eqref{equation:pump-as-load}\\
%     & & & \textnormal{\color{blue}{Budget constraints}}.
% \end{aligned}
% \end{equation}
% {\color{red} Here, I am considering an objective of maximizing the weighted sum of power and water demands with a constraint on the expansion budget.}

\section{Numerical Experiments}
\label{sec:results}
\begin{figure}
    \centering
    \includegraphics[width = 0.5\textwidth]{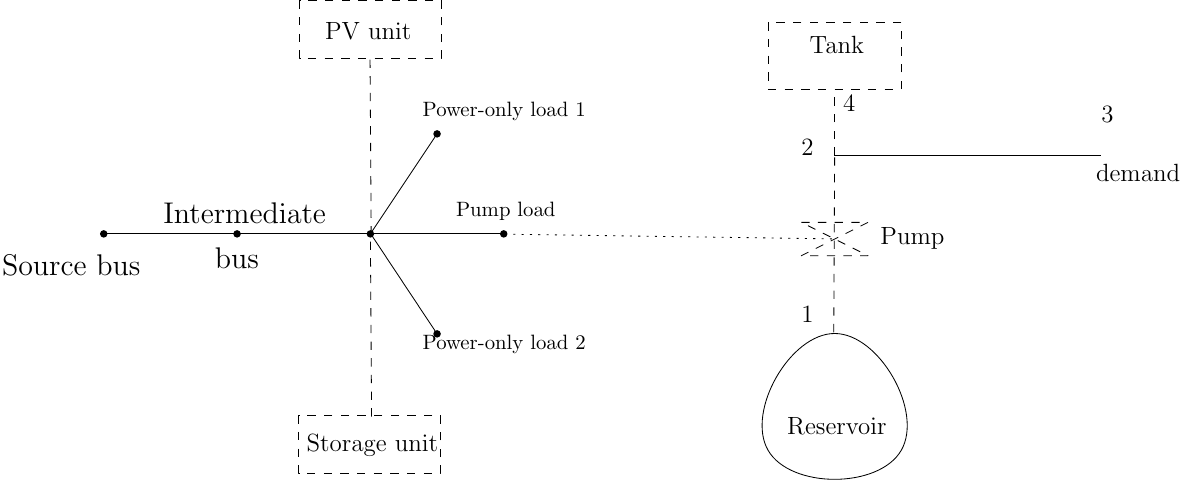}
    \caption{Illustration of integrated power (left) and water (right) networks, with dependencies (dotted lines) and expandable components (dashed lines).}
    \label{fig:joint-network}
\end{figure}

\begin{figure}
     \centering
     \begin{subfigure}[b]{0.24\textwidth}
         \centering
         \includegraphics[width=\textwidth]{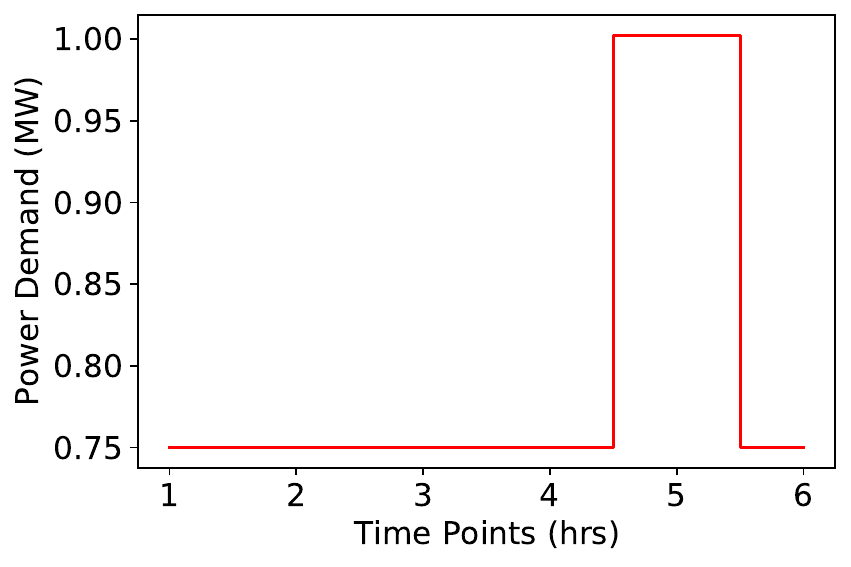}
         \caption{Power demand}
         \label{fig:power-demand}
     \end{subfigure}
     \hfill
     \begin{subfigure}[b]{0.24\textwidth}
         \centering
         \includegraphics[width=\textwidth]{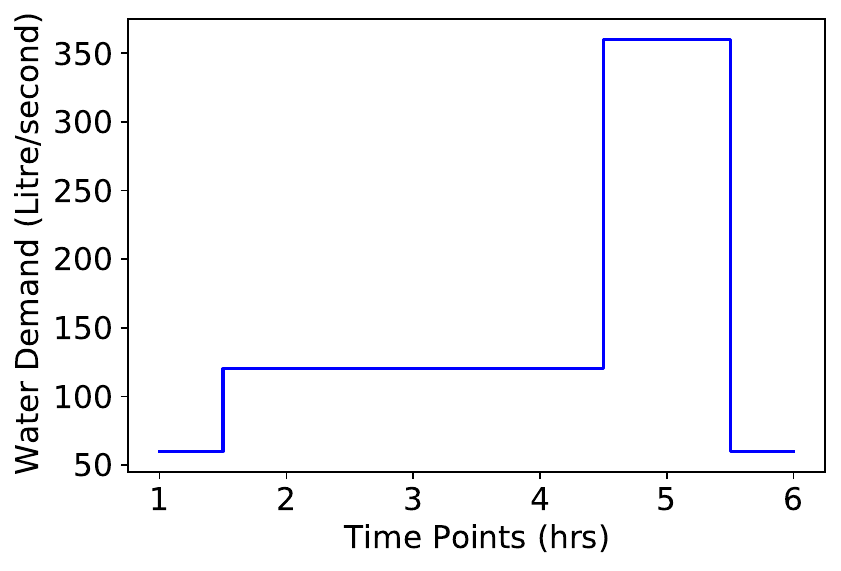}
         \caption{Water demand}
         \label{fig:water-demand}
     \end{subfigure}
    \caption{The figure illustrates the power and water demand profiles considered for experiments presented in this work. The demands peak simultaneously at hour 5.}
    \label{fig:demands}
\end{figure}

In this section, we perform numerical experiments using the aforementioned models and formulation to demonstrate the benefits of joint expansion planning. To facilitate the demonstration, we implemented the Mixed Integer Nonlinear Program (MINLP) formulation in the language Julia, using the package JuMP. We solved the formulation using SCIP \cite{scip_ref}, a non-commercial solver for MINLPs. For the experiments, we consider a joint network obtained by linking a 3-bus power distribution network with a 3-node water network, as shown in Figure \ref{fig:joint-network}. The power network consists of a source bus, a load bus (with 3 loads), and an intermediate node, with all three buses connected in a line structure. We restrict the voltage magnitudes in the network between 0.97 and 1.03 p.u. but leave the line flows unbounded. 
The water network consists of a reservoir, a demand point, and an intermediate junction. The intermediate junction and the demand point are connected by a pipe, and they are at a higher elevation compared to the reservoir. To meet the demand, 
we provide an option of adding to the network, a pump connecting the reservoir and the intermediate junction. We also provide the option of adding a tank (empty when built) at the intermediate junction.
If the pump is added, the power and water networks will be linked by connecting the new pump to one of the conductors of the load bus. For the linkage, we consider the pump curve presented in \cite{ulanicki2008modeling}.
Additionally, we provide the option of adding a PV unit and an energy storage to expand the power network.
For the demonstration, we consider a case where the power and water demand profiles are available at 1-hour intervals for a 6-hour time horizon and peak simultaneously at the $5^{th}$ hour as shown in Figures \ref{fig:power-demand} and \ref{fig:water-demand}. Note that Figure \ref{fig:power-demand} excludes pump loads. With this setup, we perform two experiments.

% {\color{blue} The distribution network used consisted of three buses connected in a line structure. The network includes a source bus, a load bus, and an intermediate node. Voltage magnitudes are bounded to be between 0.97 and 1.03 p.u., and line flows are kept unbounded.
% The water pump is connected to one of the conductors of the load bus and serves as a {\color{red}power load (talk about power-only or non-pump loads}.}
% {\color{red} Data links. Pump curve \cite{ulanicki2008modeling}, PV load? Computation times. Power-only load}

% \vspace{5pt}

% {\color{red} Add network figures, network data: pipeline, xxx etc.}

\begin{figure}
     \centering
     \begin{subfigure}[b]{0.3\textwidth}
         \centering
         \includegraphics[width=\textwidth]{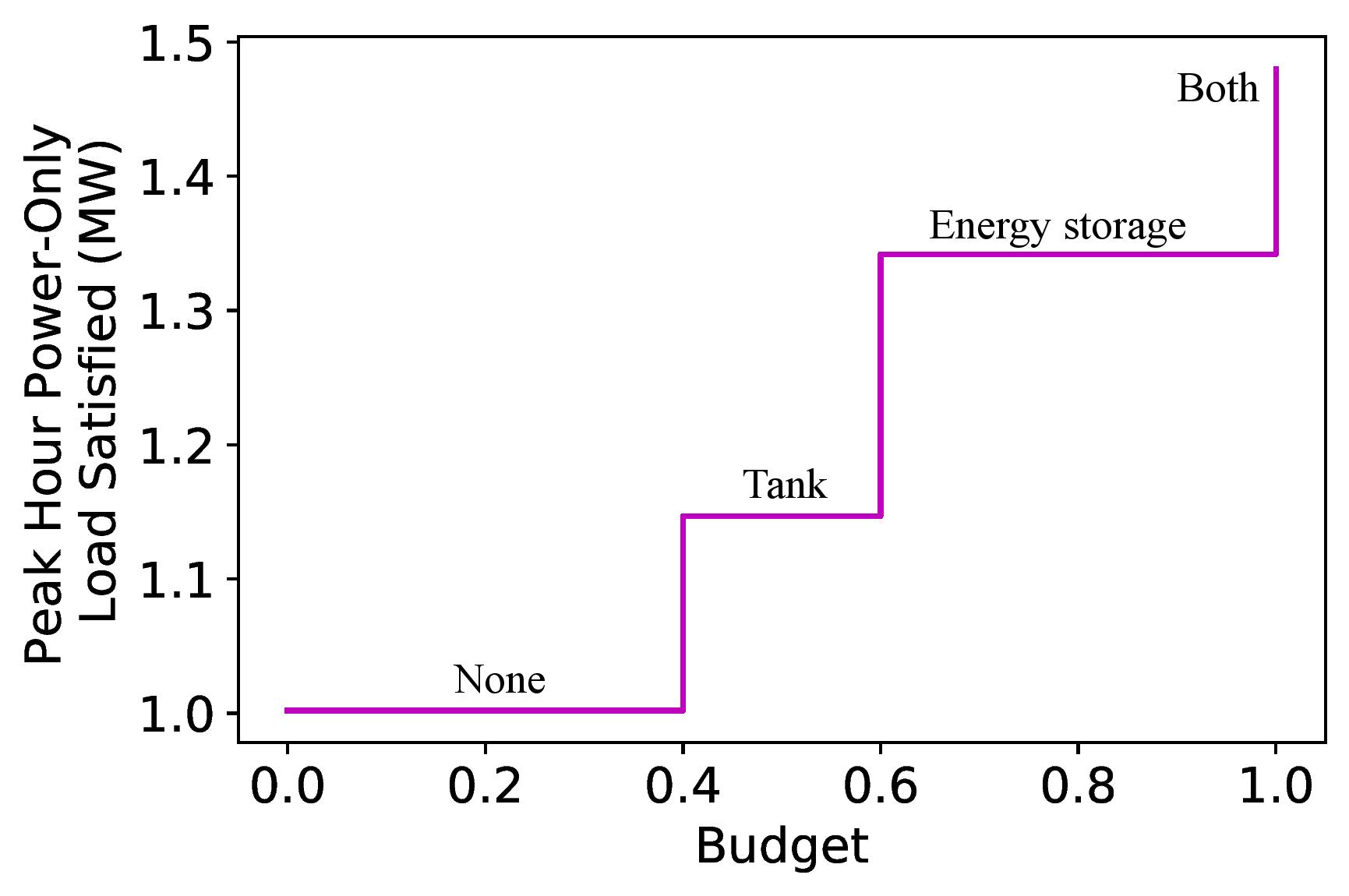}
         \caption{At budget levels 0.4, 0.6, and 1.0 respectively, a water tank, a battery storage, and both are built to meet higher power demands.}
         \label{fig:power-step-result}
     \end{subfigure}
     \hfill
     \begin{subfigure}[b]{0.3\textwidth}
         \centering
         \includegraphics[width=\textwidth]{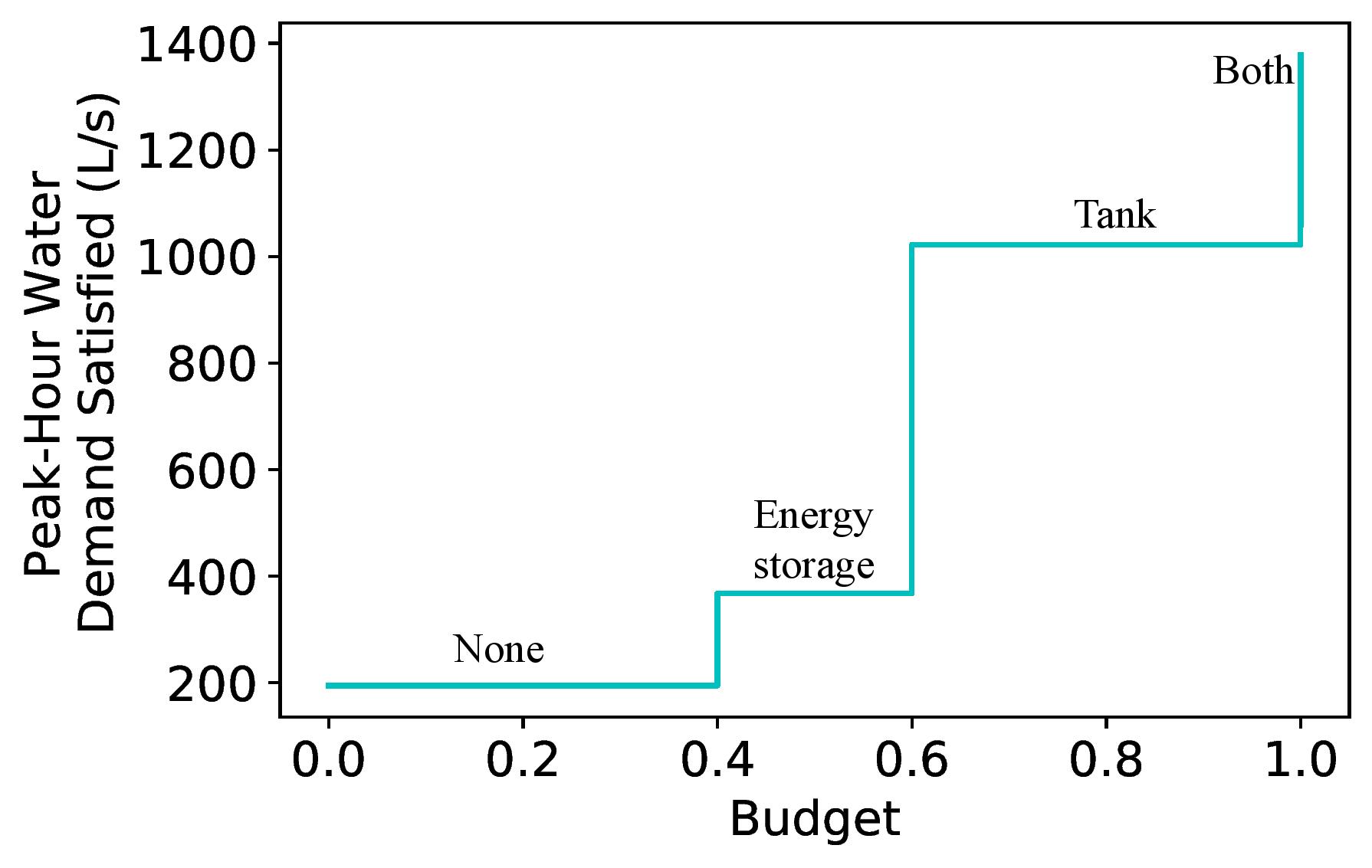}
         \caption{At budget levels 0.4, 0.6, and 1.0 respectively a battery storage, a water tank, and both are built to meet higher water demands.}
         \label{fig:water-step-result}
     \end{subfigure}
    \caption{The figure depicts the increase in the peak-hour demand met with the expansion budget.}
    \label{fig:step-budget}
\end{figure}

In the first experiment, we fix the water demand profile over all the 6 hours to the one shown in Figure \ref{fig:water-demand} and the power demand profile at all except the $5^{th}$ hour to the one described by Figure \ref{fig:power-demand}. The only variable term in the objective is the peak-hour power demand that can be met. We solve the problem of maximizing this peak-hour demand (by setting $\lambda = 1$ in the MINLP) for iteratively increasing values of the expansion budget between 0 and 1. To simulate a scenario where it is cheaper to build a tank compared to storage, we consider representative expansion costs of 0.6, 0.4, 0, and 0 units respectively for the battery storage, water tank, pump, and PV generator. The results of this experiment are summarized in Figure \ref{fig:power-step-result}.
Initially, when the budget for expansion is below 0.4, only the zero-cost components, PV generator and pump, are built to satisfy the fixed water demand and a peak power demand of 1.00 MW. When the budget is in between 0.4 and 0.6, the peak power demand is increased to 1.15 MW by building a water tank. Here, building a tank is sufficient to increase the peak-hour power demand as it helps in shifting the pumping to earlier time steps, when the power demand is low, and therefore allocates a lesser amount of power for pumping during the peak hour. When the budget reaches a value of at least 0.6, a battery storage is built to meet a higher peak demand, and at a budget of 1 unit, both storage and tank are built to satisfy the highest amount of peak power demand in the experiment.

From this experiment, it can be inferred that it is not always necessary to build an expensive component, such as a battery storage in this case, especially when a cheaper option such as a water tank is available in the interdependent neighboring network. Therefore, by offering an increased number of options to choose from, combining the expansion planning of these interdependent networks reduces the cost of expansion and potentially leads to a significant amount of savings when implemented at a large scale.

In the second experiment, we switch the costs of the battery and tank to 0.4 and 0.6 respectively. Then, we fix the peak hour power demand to 1.05 MW, which is a value in between 1.00 MW and 1.15 MW, the peak demands satisfiable with and without building a tank in the previous experiment. Instead, we let the peak-hour water demand to be a variable and maximize this peak demand by iteratively solving the MINLP formulation  (with $\lambda = 0$) with an increasing value of budget between 0 and 1. In this case, the peak water demand satisfiable is only 194.07 L/s, which is less than the fixed peak water demand satisfiable in the previous experiment (360 L/s). This is because a lower amount of power is available for pumping. The results of this experiment are summarized in Figure \ref{fig:water-step-result}, and they follow a similar trend as the previous experiment. A cheaper battery storage is built before a water tank to increase the peak-hour water demand satisfiable, even though the battery is from a neighboring interdependent network. 

Besides, suppose that the networks are expanded independently by assigning a fixed amount of power for pumping based on the requirement of other loads. Then, to satisfy an elevated demand combination of 1.05 MW and 360 L/s, independent expansion planning requires the building of two components, a battery storage to meet the peak power demand and a water tank to meet the peak water demand. However, as can be seen in Figure \ref{fig:water-step-result}, when the expansion planning is performed jointly, it is sufficient to build either one of these components to satisfy the elevated levels of demand. Therefore, the joint expansion planning not only reduces the cost of expansion, but also reduces the redundancy in expansion. In addition to cost reduction, this has implications in reducing the space required for expanding the networks.

\section{Summary}
In this paper, we formulated the joint expansion planning of power and water distribution networks as an MINLP, where the dependency between the networks arises through the power consumption of pumps. By solving the formulation on a small-scale test network, we demonstrated that an elevated power demand can be satisfied by constructing a water tank and an increased water demand can be satisfied by building an energy storage unit. Through this, we showed that the joint expansion planning offers advantages such as lower cost and reduced redundancy over expanding the networks independently. In the future, to realize the full potential of the joint planning, this formulation needs to be implemented and solved on networks of larger scale. However, the non-convexity of the problem presents computational challenges, and one needs to develop efficient algorithms to address this challenge.

\bibliographystyle{IEEEtran}
\bibliography{references.bib}

\end{document}